\documentclass[reqno]{amsart}
\usepackage{hyperref}
\usepackage{amssymb}
\usepackage{cite}

\DeclareMathOperator{\Div}{Div}
\DeclareMathOperator{\grad}{\nabla}

\DeclareMathOperator{\lip}{Lip}
\DeclareMathOperator{\op}{op}
\newcommand{\R}{\mathbb{R}}
\newcommand{\rn}{{\mathbb{R}^n}}
\newcommand{\bg}{\mathbf{g}}
\newcommand{\bu}{\mathbf{u}}

\begin{document}
\title[Bounded weak solutions]
{Bounded weak solutions with Orlicz space data: an overview}

\author[Cruz-Uribe\hfil \hfilneg] {David Cruz-Uribe, OFS}  % in alphabetical order

\address{David Cruz-Uribe, OFS\newline
Dept. of Mathematics \\
University of Alabama \\
 Tuscaloosa, AL 35487, USA}
\email{dcruzuribe@ua.edu}


\begin{abstract}
It is well known that non-negative solutions to the Dirichlet
 problem $\Delta u =f$ in a bounded domain $\Omega$, where $f\in L^q(\Omega)$, $q>\frac{n}2$, satisfy
 $\|u\|_{L^\infty(\Omega)} \leq C\|f\|_{L^q(\Omega)}$.  We generalize
 this result by replacing the Laplacian with a degenerate elliptic
 operator, and we show that we can take the data $f$ in an Orlicz space
 $L^A(\Omega)$  that, in the classical case, lies strictly between $L^{\frac{n}{2}}(\Omega)$ and
 $L^q(\Omega)$, $q>\frac{n}2$.
\end{abstract}

\maketitle \numberwithin{equation}{section}
\newtheorem{theorem}{Theorem}[section]
\newtheorem{corollary}[theorem]{Corollary}
\newtheorem{lemma}[theorem]{Lemma}
\newtheorem{remark}[theorem]{Remark}
\newtheorem{problem}[theorem]{Problem}
\newtheorem{example}[theorem]{Example}
\newtheorem{definition}[theorem]{Definition}
\allowdisplaybreaks

\section{Introduction: uniformly elliptic operators}

In this note we survey recent results from~\cite{MR4280269}, done
jointly with Scott Rodney.  Let
$\Omega$ be a bounded, open, and connected subset of $\R^n$, let
$Q=Q(x)$ be an $n\times n$ positive semi-definite, self-adjoint
measurable matrix function, and let $v$ be a non-negative, measurable
function. (Hereafter, we refer to $v$ as a weight.)  We are interested
in  studying solutions of the Dirichlet problem for the
elliptic PDE
\begin{eqnarray}\label{dp}
\left\{\begin{array}{rcll}
-\Div\left(Q\nabla u\right)&=&vf&\textrm{for }x\in\Omega\\
u&=&0&\textrm{for }x\in\partial\Omega
\end{array}\right.
\end{eqnarray}
Provided that $v(x)>0$ a.e. we can define the operator
\[ Lu = -v^{-1}\Div\left(Q\nabla u\right), \]
and while this does not make sense if $v(x)=0$ on a set of positive
measure, we will abuse notation and say that we are interested in
solutions of the equation $Lu=f$.  As we will see, there is a
close interaction between the matrix $Q$ and the weight $v$; in turn,
much of the work on this equation is informed by the theory of
weighted norm inequalities in harmonic analysis.

Our work is
motivated by earlier results by Fabes, Kenig and
Serapioni~\cite{MR643158}, Chanillo and Wheeden~\cite{MR805809},
Franchi, Lu, and Wheeden~\cite{MR1354890}, and Sawyer and
Wheeden~\cite{SW2,MR1175693,MR2204824}.   Our previous work with
Rodney and Rosta~\cite{CRR2,CRR1} is also relevant.

Our starting point is the following classical result for uniformly
elliptic operators  due to
Trudinger~\cite{MR369884} (see also
Maz$^\prime$ya~\cite{MR0131054,MR0259329} and
Stampacchia~\cite{MR192177}).

  \begin{theorem} \label{thm:trudinger}
    Let $f\in L^q(\Omega)$, $q>\frac{n}{2}$, $Q$ uniformly elliptic,
    and $v=1$.  If $u$ is a non-negative weak
    solution of~\eqref{dp},    then
    \[ \|u\|_{L^\infty(\Omega)} \leq C\|f\|_{L^q(\Omega)}. \]
  \end{theorem}

The standard proof of this result using Moser iteration, and the
classical Sobolev inequality
\[ \bigg(\int_\Omega
  |\psi(x)|^{\frac{2n}{n-2}}\,dx\bigg)^{\frac{n-2}{2n}}
  \leq C\bigg(\int_\Omega |\grad \psi(x)|^2\,dx\bigg)^{\frac{1}{2}}. \]
The lower bound for $q$ is closely related to the fact that
\[ \frac{n}{2} = \left(\frac{n}{n-2}\right)' \]
is the dual exponent of the ``gain'' in the Sobolev inequality.  The
bound on $q$ in Theorem~\ref{thm:trudinger} is sharp:  if
we take $Q$ to be the identity (so the operator becomes the
Laplacian) and $\Omega = B(0,1)$, and if we let
\[ f(x) = \frac{1}{|x|^{2}\log(e+|x|^{-1})}, \]
then $f\in L^{\frac{n}{2}}(\Omega)$ but the  solution to
$\Delta u =f$ is unbounded at the origin.

Our first result generalizes Theorem~\ref{thm:trudinger} by showing
that we can get closer to the endpoint by passing to a finer scale of
spaces.  Recall that a Young function $A : [0,\infty)\rightarrow
[0,\infty)$, is an increasing, convex function that satisfies $A(0)=0$
and $A(t)/t\rightarrow \infty$ as $t\rightarrow \infty$.  We define
the Orlicz space $L^A(\Omega)$ to consist of all measurable functions
$f$ such that
\[ \|f\|_{L^A(\Omega)} = \inf\bigg\{ \lambda > 0 :
  \int_\Omega A\bigg(\frac{|f(x)|}{\lambda}\bigg)\,dx \leq 1 \bigg\}
  < \infty. \]
If we take $A(t) = t^{\frac{n}{2}}\log(e+t)^q$, $q>0$, then for every
$\epsilon>0$, 
\[ L^{\frac{n}{2}+\epsilon}(\Omega)
  \subsetneq
  L^A(\Omega)
  \subsetneq
  L^{\frac{n}{2}}(\Omega). \]
For more information on Orlicz spaces, see~\cite{RR}.

\begin{theorem} \label{thm:orlicz-trudinger}
  Let $f\in L^A(\Omega)$, $A(t)= t^{\frac{n}{2}}\log(e+t)^q$,
  $q>\frac{n}{2}$, $Q$ uniformly elliptic, and $v=1$.
If $u$ is a non-negative weak
    solution of~\eqref{dp},    then
    \[ \|u\|_{L^\infty(\Omega)} \leq C\|f\|_{L^q(\Omega)}. \]
  \end{theorem}

  Theorem~\ref{thm:orlicz-trudinger} is a special case of our main
  result, Theorem~\ref{thm:degen-trudinger} below.  It is, however,
  not new:  earlier, Cianchi~\cite{MR1721824} proved it using very
  different techniques; he also gave a better lower bound, proving it
  for $q>\frac{n}{2}-1$.  The above example shows that this bound is
  sharp.

  \section{Degenerate elliptic operators}

  Our main result is a generalization of
  Theorem~\ref{thm:orlicz-trudinger} that holds for a large class of
  degenerate elliptic operators.  Our approach, following our previous
  work in~\cite{CRR2,CRR1} and the earlier work of Sawyer and
  Wheeden~\cite{MR2204824,SW2} is to give the broadest possible
  hypotheses on the matrix $Q$ and the weight $v$ for which our
  results hold.  We make three critical assumptions:
  \begin{itemize}

  \item $v \in L^1(\Omega)$;

    \item $|Q(x)|_{\op} =\sup\{ |Q(x)\xi| : \xi \in \rn, |\xi|=1 \}
      \leq kv(x)$;

    \item there exists $\sigma>1$ such that for all
      $\psi \in \lip_0(\Omega)$ (that is, compactly supported Lipschitz
      functions)
        \[ \bigg(\int_\Omega
  |\psi(x)|^{2\sigma}v(x)\,dx\bigg)^{\frac{1}{2\sigma}}
  \leq C\bigg(\int_\Omega |Q^{\frac{1}{2}}(x)\grad \psi(x)|^2\,dx\bigg)^{\frac{1}{2}}. \]
\end{itemize}

These hypotheses hold in a number of cases: if we take $v=1$ and let
$Q$ be uniformly elliptic, the constant $k$ is just the upper bound on
the eigenvalues of $Q$, and we can take $\sigma = \frac{n}{n-2}$ and
use the
  classical Sobolev inequality.  Fabes, Kenig and
  Serapioni~\cite{MR643158} considered the case when the weight $v$ satisfies the
  Muckenhoupt $A_2$ condition,
  \[ [v]_{A_2} = \sup_B
  \frac{1}{|B|}\int_B v(x)\,dx \;\frac{1}{|B|}\int_B v(x)^{-1}\,dx
  < \infty, \]
where $B$ is any ball, and $Q$ satisfies the degenerate
ellipticity condition
\[ \lambda v(x)|\xi|^2 \leq \langle Q(x)\xi,\xi\rangle \leq
  \Lambda  v(x)|\xi|^2. \]
They proved that in this case the Sobolev inequality holds for
$\sigma=\frac{n}{n-1}+\delta$, where $\delta>0$ depends on $n$ and
$[w]_{A_2}$.

Chanillo and Wheeden~\cite{MR805809} introduced the concept of
$2$-admissible pairs.  They considered matrices $Q$ that satisfy
\[ w(x)|\xi|^2 \leq \langle Q(x)\xi,\xi\rangle \leq
  v(x)|\xi|^2, \]
where $w,\,v$ are weights that satisfy $w(x)\leq v(x)$, $v$ doubling,
$w\in A_2$, and together satisfy a balance condition: there
exists $\sigma>1$ such that given $B_1\subset B_2 \subset \Omega$,
\begin{equation} \label{eqn:balance}
 \frac{r(B_1)}{r(B_2)}
  \bigg(\frac{v(B_1)}{v(B_2)}\bigg)^{\frac{1}{2\sigma}}
  \leq C \bigg(\frac{w(B_1)}{w(B_2)}\bigg)^{\frac{1}{2}}, 
\end{equation}
where $r(B)$ is the radius of the ball $B$.  They proved that in this
case a weighted Sobolev inequality holds:
\begin{multline*}
 \bigg(\int_\Omega
 |\psi(x)|^{2\sigma}v(x)\,dx\bigg)^{\frac{1}{2\sigma}} \\
 \leq
 C\bigg(\int_\Omega |\grad \psi(x)|^2 w(x)\,dx\bigg)^{\frac{1}{2}}
  \leq\bigg(\int_\Omega |Q^{\frac{1}{2}}(x)\grad
  \psi(x)|^2\,dx\bigg)^{\frac{1}{2}}.
\end{multline*}
Explicit examples of weights that satisfy the balance condition~\eqref{eqn:balance} can
be found in~\cite{MR2771262}.

\smallskip

The next step is to give a definition of weak solutions that is
adapted to our operator.  We follow the approach developed
in~\cite{CRR1}.   Define the degenerate Sobolev space
$QH_0^1(v;\Omega)$ to be the closure of $\lip_0(\Omega)$ with respect
to the norm
  \begin{multline*}
    \|\psi\|_{QH^1_0(v;\Omega)} = \|\psi\|_{L^2(v;\Omega)} + \|\nabla
    \psi\|_{L^2_Q(\Omega)} \\
    =
    \bigg(\int_\Omega |\psi|^2v\,dx\bigg)^{\frac{1}{2}}
    +
    \bigg(\int_\Omega |\sqrt{Q(x)}\nabla \psi|^2\,dx\bigg)^{\frac{1}{2}}.
  \end{multline*}

Formally,  $QH_0^1(v;\Omega)$ consists of equivalence classes of
Cauchy sequences; however, we can define a unique pair of functions
that is associated with each such class.  Given a Cauchy sequence
$\{\phi_k\}_{k=1}^\infty$, since both $L^2(v;\Omega)$ and
$L^2_Q(\Omega)$ are Banach spaces, the sequence converges to some
function $u$ in $L^2(v;\Omega)$, and the sequence $\{\grad
\phi_k\}_{k=1}^\infty$ converges to some vector-valued function $\bg$
in $L^2_Q(\Omega)$.  We will write $\grad u = \bg$ and think of it as
the weak gradient of $u$.  However, it is important to note that in
this setting, $\bg$ may not be a weak derivative in the classical
sense if the matrix $Q$ is too degenerate.  In~\cite{MR643158}, the
authors give an example of a matrix $Q$ and a pair $(u,\bg)$ such that
$u$ is non-constant, but $\bg=0$.

We now define weak solutions of the Dirichlet problem~\eqref{dp} to be
any pair $\bu = (u,\grad u) \in QH_0^1(v;\Omega)$  that satisfies
\[ \int_\Omega \nabla \psi(x) \cdot Q(x)\nabla u (x) \,dx
  = \int_\Omega f(x)\psi(x)v(x)\,dx \]
for every  $\psi\in \lip_0(\Omega)$.

\smallskip

We can now state our main result.  Here $L^A(v;\Omega)$ is an Orlicz
space defined as above but with respect to the measure $v\,dx$. 

\begin{theorem} \label{thm:degen-trudinger}
Let $Q$ and $v$ satisfy the above hypotheses with gain $\sigma>1$ in
the Sobolev inequality.  Let $f\in
L^A(v;\Omega)$, where $A(t)=t^{\sigma'} \log(e+t)^q$, $q>\sigma'$.  If
$\bu=(u,\grad u) \in QH_0^1(v;\Omega)$ is a non-negative weak
solution of~\eqref{dp},  then
\begin{equation*}
\|u\|_{L^\infty(v;\Omega)} \leq C\|f\|_{L^A(v;\Omega)}.
\end{equation*}
  \end{theorem}

  The proof of Theorem~\ref{thm:degen-trudinger} is loosely modeled on
  one of the proofs of Theorem~\ref{thm:trudinger}.  As we noted above, the
  typical proof of this result uses Moser iteration, but we were
  unable to make it work in our setting.   Instead, we used a proof
  that relied on De Giorgi iteration, adapting ideas from the recent
  work of Korobenko, {\em et al.}~\cite{Korobenko:2016ue}.

  The first step in the proof is technical:  as we noted above, given
  a pair $(u,\grad u) \in QH_0^1(v;\Omega)$, $\grad u$ may not be a
  classical weak derivative.  Nevertheless, we need it to satisfy many
  of the same properties; in proving that they do, we make very heavy
  use of the hypothesis that $|Q|_{\op} \leq kv$.  This should be
  contrasted with our results in~\cite{CRR1} which did not require
  this assumption.

  Given these properties, we can now begin the process of De Giorgi
  iteration.  For each $r>0$, define
\[ \phi = \phi_r(u) = (u-r)_+. \]
Let $ S(r) = \{ x : u(x)> r \}$.  Then we have that

\[ (\phi,\nabla \phi) = ((u-r)_+, \chi_{S(r)} \nabla u) \in
  QH_0^1(\Omega). \]
By the Sobolev inequality we assume to hold, and by the definition of
weak solutions (using  $\phi$ as our test function), we have that
\[ \|\phi\|_{L^{2\sigma}(v;\Omega)}^2 
\leq 
C_0\|f\|_{L^{(2\sigma)'}(v;\Omega)}\|\phi\|_{L^{2\sigma}(v;S(r))}. \]
By H\"older's inequality in the scale of Orlicz spaces,  and by the
definition of the Orlicz norm,
for $s>r$ we have that
\begin{align} \label{eqn:iteration}
 v(S(s))^{\frac{1}{2\sigma}}(s-r)
& \leq \|\phi\|_{L^{2\sigma}(v;\Omega)} \\
& \leq C\|f\|_{L^{(2\sigma)'}(v;S(r))} \notag \\
& \leq
C\|f\|_{L^A(v;\Omega)}  \|\chi_{S(r)}\|_{L^B(v;\Omega)} \notag \\
& \leq C \|f\|_{L^A(v;\Omega)}
 \frac{v(S(r))^{\frac{1}{2\sigma}}}
{\log(e+(v(S(r)))^{-1})^{q\left(\frac{(2\sigma)'}{\sigma'}\right)}}  \notag
\end{align}

If we now define
\[ C_k = \tau_0\|f\|_{L^A(v;\Omega)}\bigg(1-
  \frac{1}{(k+1)^\epsilon}\bigg),  \]
then to complete the proof we need to show that
\[ v( S(\tau_0\|f\|_{L^A(v;\Omega)})) 
= 
\lim_{k\rightarrow\infty} v(S(C_k))
=
0. \]
To prove this, let $\displaystyle m_k = -\log(v(S(C_k))$; then the
above limit is equivalent to showing that $m_k\rightarrow \infty$ as
$k\rightarrow \infty$.

In inequality~\eqref{eqn:iteration}, let $s=C_{k+1}$, $r=C_k$.  If we fix
$\epsilon= \frac{q}{\sigma'}-1>0$, then
\[ m_{k+1} \geq
  \log\left(\displaystyle\frac{\epsilon\tau_0}{C}\right)
  + \log\left(\displaystyle\frac{m_k}{k+2}\right)^\frac{2\sigma
    q}{\sigma'}
  +m_k. \]
By induction, we can show that there exists $\tau_0>0$ (very large) such that
$m_k \geq m_0 + k$.
  Therefore,
  \[ \lim_{k\rightarrow \infty} m_k = \infty. \]

  \section{Further remarks}

 Theorem~\ref{thm:degen-trudinger} is part of a larger project to
 develop a theory of existence, uniqueness, and regularity for
 degenerate PDEs.  We want to close this note by outlining some
 further directions.  Motivated directly by the work
 in~\cite{MR4280269}, we see three immediate problems.  First, as we
 noted, Theorem~\ref{thm:orlicz-trudinger}, which is just a special
 case of Theorem~\ref{thm:degen-trudinger}, is not sharp.  We have
 not, despite repeated efforts, been able to improve our
 argument.  Therefore, it is open whether
 Theorem~\ref{thm:degen-trudinger} can be improved so that in the
 classical case we get a sharp result.  Second, and perhaps related,
 our proof uses De Giorgi iteration.  We originally attempted to use
 Moser iteration but were not successful.  It would be interesting to
 make Moser iteration work in this setting.  Towards this end, we have
 generalized the classical identiy
 \[ \lim_{p\rightarrow \infty} \|f\|_{L^p(\Omega)} =
   \|f\|_{L^\infty(\Omega)}, \]
 which plays a central role in Moser iteration, to the scale of Orlicz
 spaces:  see~\cite{MR4556453}.   Third, the classical result is known
 for equations with lower order terms, so it should be possible to
 formulate and prove a similar result for degenerate equations.

 Looking beyond this, there are several directions we believe are
 worth exploring.  The first is whether we can extend these results to
 data in other function spaces, such as Lorentz spaces, the small
 Lebesgue spaces (that is, the dual spaces of the grand Lebesgue
 spaces of Iwaniecz and Sbordone--see~\cite{MR2110048}), or the
 variable Lebesgue spaces.  Each of these families of function spaces
 would provide new insight into the behavior of solutions as the data
 function gets close to the endpoint space.

 Second, we want to examine whether the hypothesis that we have a
 Sobolev inequality with gain $\sigma$ can be weakened.  Motivated by
 problems in the study of hypoelliptic operators, Korobenko, {\em et
   al.}~\cite{Korobenko:2016ue} have introduced Sobolev inequalities
 with gain in the scale of Orlicz spaces:  roughly, replacing the
 $L^{2\sigma}$ norm on the lefthand side with an Orlicz norm given by
 $A(t)=t^2\log(e+t)^\sigma$.   We think that it is worth exploring
 whether Theorem~\ref{thm:degen-trudinger} can be extended in this
 direction.

 Third, having given conditions when solutions are bounded, the next
 step is to see whether we can prove regularity of solutions.  At the
 heart of this problem will be to determine whether additional
 hypotheses are required.  Preliminary results in this direction were
 obtained in~\cite{MR2771262,MR3011287}.  

Finally, we are interested in understanding further the hypothesis
that a Sobolev inequality (either in the scale of Lebesgue or Orlicz
spaces) exists.  Given a matrix $Q$ and a weight $v$, what are
necessary and/or sufficient conditions on the domain $\Omega$ for some
kind of Sobolev inequality to exist?  We are
particularly interested in generalizing the geometric
characterizations in~\cite{MR2386936} (see also~\cite{MR4107806}).

\section*{Acknowledgements}
The author acknowledges support from research funds from the Dean of the
College of Arts \& Sciences, the University of Alabama, and from a
Simons Foundation Travel Support for Mathematicians  Grant.

\bibliographystyle{plain}
\bibliography{Ghent}

\end{document}